\documentclass[12 pt]{amsart}
\usepackage{verbatim}
\usepackage{amssymb}

\usepackage{enumerate}
\usepackage[active]{srcltx} 
\numberwithin{equation}{section}
\usepackage{hyperref}

\usepackage{t1enc}
\usepackage[utf8x]{inputenc}

\newtheorem{theorem}{Theorem}[section]

\newtheorem*{theorem*}{Theorem}%[section]

\newtheorem{lemma}[theorem]{Lemma}
\newtheorem{problem}[theorem]{Problem}
\newtheorem{proposition}[theorem]{Proposition}

\newtheorem{corollary}[theorem]{Corollary}

\newcommand{\wreso}[1]{monotonically $#1$-resolvable}

\theoremstyle{definition}
\newtheorem{definition}[theorem]{Definition}
\theoremstyle{remark}

\newcommand{\mbb}[1]{\mathbb{#1}}

\newcommand{\setm}{\setminus}

\newcommand{\subs}{\subset}

\def\<{\left\langle}
\def\>{\right\rangle}

\author[I. Juh\'asz]{Istv\'an Juh\'asz}
\thanks
  {
   }
\address
      { Alfr{\'e}d R{\'e}nyi Institute of Mathematics , Eötvös Loránd Research Network
}
\email{juhasz@renyi.hu}

\author[L. Soukup]{Lajos Soukup}
\thanks
  {
   }
\address
      { Alfr{\'e}d R{\'e}nyi Institute of Mathematics, Eötvös Loránd Research Network
}
\email{soukup@renyi.hu}

\author[Z. Szentmikl\'ossy]{Zolt\'an Szentmikl\'ossy}
\address{Eötvös University of Budapest}
\email{szentmiklossyz@gmail.com}


\title{On resolvability of products}
\thanks{The research on and preparation of this paper was
supported by  NKFIH grant no. K129211}
\date{\today}

\begin{document}

\begin{abstract}
All spaces below are $T_0$ and crowded  (i.e. have no isolated points).

For $n \le \omega$ let $M(n)$  be the statement that there are $n$ measurable cardinals
and $\Pi(n)$ ($\Pi^+(n)$) that there are $n+1$ (0-dimensional $T_2$) spaces whose product is irresolvable.
We prove that $M(1),\,\Pi(1)$ and $\Pi^+(1)$ are equiconsistent. For $1 < n < \omega$ we show that
$CON(M(n))$ implies $CON(\Pi^+(n))$. Finally, $CON(M(\omega))$ implies the consistency of having
infinitely many crowded 0-dimensional $T_2$-spaces such that the product of any finitely many of them is irresolvable.
These settle old problems of Malychin from \cite{Ma73}.

Concerning an even older question of Ceder and Pearson in \cite{CePe},  we show that the following
are consistent modulo a measurable cardinal:

\begin{enumerate}[(i)]
\item There is a 0-dimensional $T_2$ space $X$ with $\omega_2 \le \Delta(X) \le 2^{\omega_1}$
whose product with any countable space is not $\omega_2$-resolvable, hence not maximally resolvable.

\smallskip

\item There is  a monotonically normal space $X$ with $\Delta(X) = \aleph_\omega$
whose product with any countable space is not $\omega_1$-resolvable, hence not maximally resolvable.

\end{enumerate}
These significantly improve a result of Eckertson in \cite{Eck97}.

\end{abstract}

\maketitle

\section{Introduction}

It is an easy exercise to show that any product of infinitely many non-singleton spaces
is $\mathfrak{c}$-resolvable. On the other hand, the question if there are two
crowded spaces whose product is irresolvable turned out to be non-trivial. Malychin noticed
in \cite{Ma73} that if $\mathcal{U}$ is a countably complete ultrafilter on $\kappa$ then
the product of the, obviously crowded, $T_1$-space $\<\kappa, \mathcal{U} \cup \{\emptyset\}\>$
with any countable irresolvable space is irresolvable. So, the existence of a measurable cardinal
yields an affirmative answer to this question.

This result naturally lead him to ask the following
questions:

\begin{enumerate}[(M1)]
  \item  Can $T_1$ be improved to $T_2$ or $T_3$?

\smallskip

  \item  Can we find three crowded spaces whose product is irresolvable?
\end{enumerate}

In the first part of this paper we answer both questions affirmatively by proving the following results.

\begin{enumerate}

\item The existence of a measurable cardinal is equiconsistent with having two crowded 0-dimensional
  $T_2$-spaces whose product is irresolvable.

\smallskip

\item For each $n \in \omega$, the existence of $n$ measurable cardinals implies the consistency of the existence of $n+1$ crowded 0-dimensional
  $T_2$-spaces whose product is irresolvable.
(We do not know if equiconsistency holds for $n > 1$.)

\smallskip

\item The existence of infinitely many measurable cardinals implies that it is consistent to have infinitely many crowded 0-dimensional
$T_2$-spaces such that the product of any finitely many of them is irresolvable.
\end{enumerate}

\medskip

In the second part we present some new results concerning the following old question of Ceder and Pearson:
Is the product of a maximally resolvable space with any other space maximally resolvable? By improving results of Eckertson from \cite{Eck97}
we show that the following two statements are consistent, again modulo a measurable cardinal:
\begin{enumerate}[(i)]
\item There is a 0-dimensional $T_2$ space $X$ with $\omega_2 \le |X| = \Delta(X) \le 2^{\omega_1}$
whose product with any countable space is not $\omega_2$-resolvable.

\smallskip

\item There is  a monotonically normal space $X$ with $|X| = \Delta(X) = \aleph_\omega$
whose product with any countable space is not $\omega_1$-resolvable.
\end{enumerate}
Since countable $\omega$-resolvable, hence maximally resolvable spaces exist, both results yield counterexamples
to the question of Ceder and Pearson.

We do not know if the assumption of the existence of a measurable cardinal is really needed in
these results or, for that matter, in finding any counterexample to the question of Ceder and Pearson.

\section{Irresolvable products}

To save characters, in the rest of this paper "space" will always stand for "crowded $T_0$-space".
In particular, this implies that any non-empty open set in such a space is infinite.

Since the subject of this section is "irresolvable products", we start by pointing out that
such a product has to have only finitely many factors.
Since, by our above convention, no space is a singleton, this is immediate from the following proposition.

\begin{proposition}\label{pr:inf}
If a product space has infinitely many factors then it is $\mathfrak{c}$-resolvable.
\end{proposition}

\begin{proof}
Clearly, it suffices to show that any product of the form $X = \prod\{X_n : n < \omega\}$ is $\mathfrak{c}$-resolvable.
Consider the equivalence relation $\sim$ on $X$ for which $x \sim y$ iff $\{n : x(n) \ne y(n)\}$ is finite.
Then the equivalence class $[x]$ of any point $x \in X$ is dense in $X$, moreover it is obvious that
there are at least $\mathfrak{c}$ distinct, hence disjoint equivalence classes.
\end{proof}

Next, we are going to give a condition implying that the product of two spaces is resolvable. So, its
failure is a useful necessary condition for the product of two spaces to be irresolvable.
For this we need to introduce the following notion that is clearly a weakening of $\kappa$-resolvability.

\begin{definition}
For an infinite cardinal $\kappa$, the monotone increasing family $\{A_\alpha : {\alpha}<{\kappa}\}$
is said to be a {\em monotone $\kappa$-resolution} of
the space $X$ if $$X = \bigcup \{A_\alpha : {\alpha}<{\kappa}\} \text{ and  }
int(A_\alpha) = \emptyset \text{ for each }{\alpha}<{\kappa}.$$
In this case we also say that $X$ is {\em monotonically $\kappa$-resolvable}.
\end{definition}
This concept was introduced in \cite{SoSt}, using a different but equivalent definition.
It was used there in the study of resolvability properties of spaces from the ground model
in generic extensions. But, as we shall see, it turns out to be quite useful in the study of resolvability
properties of products as well.

It is useful to note that $\{A_\alpha : {\alpha}<{\kappa}\}$ is a {\em monotone $\kappa$-resolution} of $X$ iff
the sequence of complements $\{X \setm A_\alpha : {\alpha}<{\kappa}\}$ is a {\em decreasing} $\kappa$-sequence
of sets {\em dense} in $X$ with empty intersection.

We are going to use the following piece of notation for any space $X$:
$$MR(X) = \{\kappa : X \text{ is \wreso{{\kappa}}} \}.$$
Of course, $MR(X)$ can be empty, however if $X$ is neat, i.e. $|X| = \Delta(X)$, then $|X| \in MR(X)$.
Also, if $X$ is of first category (in itself) then $\omega \in MR(X)$.

\begin{theorem}\label{tm:MRres}
If $MR(X) \cap MR(Y) \ne \emptyset$ then $X \times Y$ is resolvable.
\end{theorem}

\begin{proof}
Assume that $\kappa \in MR(X) \cap MR(Y)$ and $X = \bigcup \{A_\alpha : {\alpha}<{\kappa}\}$, $Y = \bigcup \{B_\alpha : {\alpha}<{\kappa}\}$ are
witnessing this. For each point $x \in X$ let $$rk_X(x) = \min\{\alpha : x \in A_\alpha\},$$
and define $rk_Y(y)$ for $y \in Y$ analogously.

Now, let us define the 2-coloring $c : X \times Y \to 2$ by $c(x,y) = 0$ iff $rk_X(x) \le rk_Y(y)$.
Since every non-empty open set in $X$ (resp. $Y$) meets $A_\alpha$ (resp. $B_\alpha$) for cofinally many $\alpha$
below $\kappa$, it is obvious that both $c^{-1}(0)$ and $c^{-1}(1)$ are dense in $X \times Y$.
\end{proof}

It is immediate from this that if $X$ and $Y$ are neat spaces with $|X| = |Y|$ then then $X \times Y$
is resolvable. In particular, $X^2$ is resolvable whenever $X$ is neat. However, we do not know the answer to
the following question.

\begin{problem}\label{pb:3res}
Is there a neat space whose square is not 3-resolvable?
\end{problem}

It also follows from Theorem \ref{tm:MRres} that if $X \times Y$ is irresolvable then $\omega \notin MR(X)$ or $\omega \notin MR(Y)$.
This simple observation is crucial in our promised claim that if an irresolvable product $X \times Y$ exists
then it is consistent to have a measurable cardinal, in fact then
there is an inner model with a measurable cardinal. To see this, we need the following result.

\begin{theorem}\label{tm:2cat}
If the space $X$ is of second category, i.e. not the union of countably many nowhere dense sets,
it has a non-empty regular open subset that, as a subspace, is Baire.
\end{theorem}

\begin{proof}
Let $\mathcal{W}$ be a maximal disjoint collection of first category open sets in $X$ and
let $W = \cup \mathcal{W}$. Clearly then $W$ is of first category and so is its closure, the regular closed
set $\overline{W}$. But then, by the maximality of $\mathcal{W}$, the non-empty regular open set
$X \setm \overline{W}$ is indeed Baire.
\end{proof}

It is important to point out that the above subspace $\overline{W}$ of $X$, if non-empty, is crowded, being regular closed in $X$.

We recall next that a space is {\em open hereditary irresolvable} (in short: OHI) if all its non-empty open subsets are irresolvable.
It is well-known that any irresolvable space has an open OHI subspace.

\begin{theorem}\label{tm:bair}
If the product $X \times Y$ is irresolvable then $X$ or $Y$ has an irresolvable Baire (regular open) subspace.
Consequently, there is an inner model with a measurable cardinal.
\end{theorem}

\begin{proof}
Any irresolvable space has an open OHI subspace, hence we have open $U \subs X$ and $V \subs Y$ such that $U \times V$
is OHI. This, in turn, implies that both $U$ and $V$ are OHI as well. But then, by Theorem \ref{tm:MRres}, we must
have $\omega \notin MR(U)$ or $\omega \notin MR(V)$; by symmetry, we may assume that $\omega \notin MR(U)$, consequently
$U$ is of second category. By Theorem \ref{tm:2cat} then $U$ has an open subset $U_0$ that is Baire. But as $U$ is OHI,
$U_0$ is also irresolvable.

Now, we only have to refer to the fact that the existence of an irresolvable Baire space implies that
there is an inner model with a measurable cardinal, see \cite{KST} and \cite{KT}.
\end{proof}

Our next result yields a sufficient condition for obtaining two spaces having irresolvable product.
This, in turn, will be used to find such a product from a measurable cardinal. For this, we first give
a definition, a piece of terminology borrowed from forcing.

\begin{definition}\label{df:lcl}
Given an uncountable cardinal $\lambda$, the family of sets $\mathcal{S}$ is called {\em $\lambda$-closed}
if, for any (limit) ordinal $\alpha < \lambda$ and {\em monotone decreasing} subfamily $\{S_\beta : \beta < \alpha\} \subs \mathcal{S}$
there is $S \in \mathcal{S}$ such that  $S \subs \bigcap \{S_\beta : \beta < \alpha\}$.
\end{definition}

\begin{theorem}\label{tm:lcl}
If the OHI space $X$ has a $\lambda$-closed $\pi$-base $\mathcal{B}$ for some $\lambda > \omega$ then the
product $X \times Y$ is irresolvable  whenever $Y$ is irresolvable with $|Y| < \lambda$.
\end{theorem}

\begin{proof}
Take $X$ and $Y$ as required, and consider any 2-partition $X \times Y = Z_0 \cup Z_1$ of their product.
We shall show that either $Z_0$ or $Z_1$ has non-empty interior in $X \times Y$.

For each $y \in Y$ then we get a 2-partition $X = Z_{y,0} \cup Z_{y,1}$ of $X$ by putting, for $i < 2$,
$$Z_{y,i} = \{x \in X : \<x, y\> \in Z_i\}.$$ Since $X$ is OHI we have then that $$int(Z_{y,0}) \cup int(Z_{y,1})$$
is dense open in $X$.

Now, fix an indexing $Y = \{y_\alpha : \alpha < \mu\}$, then by straightforward transfinite recursion,
using that $X$ is OHI and $\mathcal{B}$ is $\lambda$-closed, we may define
$B_\alpha \in \mathcal{B}$ and $i_\alpha < 2$ for all $\alpha < \mu$ such that
$\{B_\alpha : \alpha < \mu\} \subs \mathcal{B}$ is decreasing and $B_\alpha \subs int(Z_{y,i_\alpha})$.

Now, $\mu < \lambda$ implies that there is $B \in \mathcal{B}$ with $$B \subs \bigcap \{B_\alpha : \alpha < \mu\}.$$
Also, since $Y$ is irresolvable, if $A_i = \{y_\alpha : i_\alpha = i\}$ for $i < 2$ then $int(A_i) = V$
is non-empty for some $i$. But then, as $i = i_\alpha$ holds for any $y = y_\alpha \in V$,
we have $B \times V \subs Z_i$, completing our proof.
\end{proof}

Note that the $T_1$-space Malychin got from a $\sigma$-complete, hence $\omega_1$-closed, free ultrafilter,
(that we referred to in the introduction) is OHI and trivially has an $\omega_1$-closed $\pi$-base. Thus
Malychins result is an immediate consequence of Theorem  \ref{tm:lcl}.

But to find a $T_2$ or $T_3$ space with irresolvable product, as was asked by Malychin,
we need to obtain a $T_2$ or $T_3$ OHI space having a $\lambda$-closed $\pi$-base for some $\lambda > \omega$.
Actually, since any irresolvable space has an open OHI subspace, it suffices to find an irresolvable $T_2$ or $T_3$ space that has
such a $\pi$-base. Now, maximal $\lambda$-independent families yield even 0-dimensional $T_2$, hence $T_3$ such spaces, and
luckily for us, Kunen in \cite{Ku} had already constructed such families from a measurable cardinal.

\begin{definition}\label{df:mind}
A family $\mathcal{S}$ of subsets of $\kappa$ is {\em $\lambda$-independent} if for any
$\mathcal{S}_0 \in [\mathcal{S}]^{< \lambda}$ and $\mathcal{S}_1 \in [\mathcal{S}]^{< \lambda}$
with $\mathcal{S}_0 \cap \mathcal{S}_1 = \emptyset$ we have
$$\bigcap \{S : S \in \mathcal{S}_1\} \cap \bigcap \{\kappa \setm S : S \in \mathcal{S}_0\} \ne \emptyset.$$
$\mathcal{S}$ is {\em maximal $\lambda$-independent} if no proper extension of it is $\lambda$-independent.
We are going to write $M(\kappa,\lambda)$ to denote the statement that there is a
maximal $\lambda$-independent family $\mathcal{S}$ on $\kappa$ with $|\mathcal{S}| \ge \lambda$.
\end{definition}

Now, assume that $M(\kappa,\lambda)$ holds and $\mathcal{S}$ is a maximal $\lambda$-independent family on $\kappa$.
We say that $\mathcal{S}$ is {\em separating}, if for any $\{\alpha, \beta\} \in [\kappa]^2$ there is $S \in \mathcal{S}$
with $|S \cap \{\alpha, \beta\}| = 1$. By \cite{Ku}, we may assume that $\mathcal{S}$ is separating,
provided that $\kappa$ is minimal such that $M(\kappa,\lambda)$ holds for a fixed $\lambda$.

Using the notation of Kunen's book \cite{Kb}, for each partial function $p \in Fn(\mathcal{S}, 2; \lambda)$ we let
$$B_p = \bigcap \{S : p(S) = 1\} \cap \bigcap \{\kappa \setm S : p(S) = 0\},$$ and let
$$\mathcal{B}(\mathcal{S}) = \{B_p : p \in Fn(\mathcal{S}, 2; \lambda)\}.$$ Then $\mathcal{B}(\mathcal{S})$ is the
base of a 0-dimensional topology $\tau(\mathcal{S})$ that is $T_2$
because $\mathcal{S}$ is separating. Also it is obvious that $\mathcal{B}(\mathcal{S})$ is $cf(\lambda)$-closed.

Let us observe next that the topology $\tau(\mathcal{S})$ is irresolvable, this is an immediate consequence
of the maximality of $\mathcal{S}$. Indeed, for any partition $\kappa = A \cup B$ of $\kappa$, if
$\{A, B\} \cap \mathcal{S} \ne \emptyset$, say $A \in \mathcal{S}$ then clearly $B$ is not $\tau(\mathcal{S})$-dense.
But if
$\{A, B\} \cap \mathcal{S} = \emptyset$ then by maximality there is $p \in Fn(\mathcal{S}, 2; \lambda)$ such that
$B_p \subs A$ or $B_p \subs B$, hence again $A$ and $B$ cannot both be $\tau(\mathcal{S})$-dense.

Now we are ready to present the equiconsistency result from the introduction.

\begin{theorem}\label{tm:econ1}
The following statements are equiconsistent.

\begin{enumerate}
  \item There is a measurable cardinal.
  \item There is a 0-dimensional $T_2$-space whose product with any countable irresolvable
  space is irresolvable.
  \item There are two 0-dimensional $T_2$-spaces whose product is irresolvable.
  \item There are two spaces whose product is irresolvable.
\end{enumerate}
\end{theorem}

\begin{proof}
To see Con(1) $\Rightarrow$ Con(2), we first refer to Theorems 1 and 2 of \cite{Ku} which say that (1) implies the consistency
of $M(\kappa, \omega_1)$, moreover $\omega_1 < \kappa \le 2^{\omega_1}$ for the minimal such $\kappa$.
This, in turn, implies that there
is a separating maximal $\omega_1$-independent family $\mathcal{S}$ on $\kappa$. But the irresolvable 0-dimensional $T_2$ topology $\tau(\mathcal{S})$
even has an $\omega_1$-closed base, namely $\mathcal{B}(\mathcal{S})$. Consequently, Theorem \ref{tm:lcl} can be applied to
an open OHI subspace $X$ of $\<\kappa, \tau(\mathcal{S})\>$ with $\lambda = \omega_1$ to obtain (2).

Finally, (2) $\Rightarrow$ (3) and  (3) $\Rightarrow$ (4) are trivial, while Con(4) $\Rightarrow$ Con(1) immediately follows from Theorem \ref{tm:bair}.
\end{proof}

Actually, Unsolved Problem 1. of \cite{Ma73} asked for ($T_2$ or $T_3$) spaces $X$ and $Y$
of nonmeasurable cardinality (= less than the first measurable)
such that $X \times Y$ is irresolvable. Since in our above proof we have $|X|\le \kappa \le 2^{\omega_1}$ and
$Y$ countable, this requirement is satisfied.

The above $\kappa \le 2^{\omega_1}$ is, by Theorem 1 of \cite{Ku}, in some sense still very large, because it carries
a $\kappa$-complete $\omega_2$-saturated ideal.
Hence the following natural question arises that we don't know the answer to.

\begin{problem}\label{pb:size}
What is the smallest possible value of $|X \cup Y|$ for spaces $X, Y$ with $X \times Y$ is irresolvable?
\end{problem}

Let us now move from the question of irresolvability of two spaces to that of three, four, etc., finitely many spaces.
Perhaps not surprisingly, we shall get the consistency of such examples from two, three, etc., infinitely many measurable cardinals.
Actually, what we shall need is to have appropriate maximal $\lambda$-independent families on $\kappa$ for several
distinct $\lambda$ and $\kappa$ simultaneously. The following proposition that is an immediate consequence of
Theorem \ref{tm:lcl} expresses this.

\begin{proposition}\label{pr:npr}
Assume that we have cardinals $$\lambda_0 = \omega < \lambda_1 \le \kappa_1 < \lambda_2 \le \kappa_2 < ... < \lambda_n \le \kappa_n$$
for some $n < \omega$ such that $M(\kappa_i,\lambda_i)$ holds for each $0 < i \le n$. Then there are $n+1\,$
0-dimensional $T_2$-spaces whose product is irresolvable.
\end{proposition}

\begin{proof}
We can assume that  $\kappa_i$ is minimal such that $M(\kappa_i,\lambda_i)$ holds for each $0<i\le n$.
 
Let $X_0$ be any countable irresolvable 0-dimensional $T_2$-space, moreover
$\mathcal{S}_i$ be a separating maximal $\lambda_i$-independent family on $\kappa_i$ and $X_i$
be any open OHI subspace of $\<\kappa_i, \tau(\mathcal{S}_i)\>$ for each $0 < i \le n$. Then the product
$\prod\{X_i : i \le n\}$ can be shown to be irresolvable by simply applying Theorem \ref{tm:lcl} $\,n$ times.
\end{proof}

We are going to use this proposition with the choices $\lambda_i = \kappa_i$ and to do that we shall
again make use of \cite{Ku}. Namely, in the penultimate paragraph of this paper the following result is explained that
gives a recipe for obtaining ZFC models in which $M(\kappa, \kappa)$ holds, by starting with a
measurable cardinal $\kappa$.

\begin{theorem}\label{tm:Mkk}(Kunen)
Assume that GCH holds and  $\kappa$ is a measurable cardinal with $\mathcal{U}$ a normal ultrafilter on $\kappa$.
For any set $A \in \mathcal{U}$ of inaccessibles consider the reverse Easton forcing $P_A$ that adds $\lambda^+$ many
generic subsets of $\lambda$ for each $\lambda \in A$. Then $$V^{P_A} \models GCH \text{ and } M(\kappa, \kappa).$$
\end{theorem}

Let us now consider the particular case of this result where a fixed cardinal $\mu < \kappa$ is given and
$A = \{\lambda \in \mathcal{U} : \lambda \text{ is inaccessible and } \lambda \ge \mu\}.$ Then $A \in \mathcal{U}$ and
it is standard to check that $P_A$ is $\mu$-closed, hence we have $V_\mu = (V^{P_A})_\mu$.
This clearly implies that if $\nu < \mu$ is any cardinal such that $M(\nu, \nu)$
holds in $V$ then $M(\nu, \nu)$ remains valid in $V^{P_A}$. Moreover, we
also have $P_A \subs V_\kappa$ and $|P_A| = \kappa$.

This leads us to the following result that we aimed for.

\begin{theorem}\label{tm:nmes}
\begin{enumerate}[(i)]
\item Assume GCH and  $\kappa_1 < ... < \kappa_n$ be measurable cardinals. Then there is a forcing
$Q_n$ with $Q_n \subs V_{\kappa_n}$ and $|Q_n| = \kappa_n$ such that in the generic extension $V^{Q_n}$ GCH holds,
all the $\kappa_i$'s are preserved as distinct uncountable cardinals, and $M(\kappa_i, \kappa_i)$ holds for each $0 < i \le n$.

\smallskip

\item Assume GCH and that $\<\kappa_i : 1 \le i < \omega\>$ be an increasing $\omega$-sequence of measurable cardinals.
Then there is a generic extension
in which all the $\kappa_i$'s are preserved as distinct uncountable cardinals and $M(\kappa_i, \kappa_i)$ holds for each $0 < i < \omega$.
\end{enumerate}
\end{theorem}

\begin{proof}
(i) We do induction on $n$. For $n = 1$ by Kunen's theorem \ref{tm:Mkk} we just can put $Q_1 = P_{A_1}$
where $\mathcal{U}_1$ is a normal ultrafilter on $\kappa_1$ and $A_1 = \{\lambda \in \mathcal{U}_1 : \lambda \text{ is inaccessible}\}$.

Now, assume that $n > 1$ and apply the inductive hypothesis to $\kappa_1 < ... < \kappa_{n-1}$ to obtain the appropriate $Q_{n-1}$.
Then $|Q_{n-1}| = \kappa_{n-1} < \kappa_n$ implies that $\kappa_n$ remains measurable in $V^{Q_{n-1}}$, moreover we may
apply again Kunen's theorem \ref{tm:Mkk}, now in $V^{Q_{n-1}}$, to obtain the forcing $P_{A_n}$
where $\mathcal{U}_n$ is a normal ultrafilter on $\kappa_n$ and
$$A_n = \{\lambda \in \mathcal{U}_n : \lambda \text{ is inaccessible and } \lambda > \kappa_{n-1}\}.$$
Clearly, then the two-step iterated forcing $Q_n = Q_{n-1}*P_{A_n}$, which from the point of view of our original
ground model $V$ is the $n$-step iterated forcing $Q_n = P_{A_1}* ... *P_{A_n}$, is as required.

\smallskip

(ii) We may now apply the above proof of part (i) to obtain the forcing $P_{A_n}$ in $V^{Q_{n-1}}$ for each $n > 0$.
(We let $P_{A_0} = Q_0$ be the trivial forcing.) We claim that the full support iterated limit $Q$ of the sequence $\<P_{A_n} : n < \omega\>$
is as required, i.e. $V^Q \models M(\kappa_n, \kappa_n)$ for all $n > 0$.

To see this, it suffices to show that for any $k > 0$ the tail iteration of $\<P_{A_n} : k < n < \omega\>$ is $(\kappa_k)^+$-closed.
This, however is an immediate consequence of Lemma 78. from \cite{je74} because every $P_{A_n}$ is clearly $(\kappa_k)^+$-closed
in $V^{Q_{n-1}}$ whenever $k < n$.
\end{proof}

Jensen proved in \cite{je74} that any ground model $V$ has a generic extension $W$ in which GCH
holds and any measurable cardinal in $V$ will remain measurable in $W$.
Now, this result, together with Theorem \ref{tm:nmes} and Proposition \ref{pr:npr}, immediately yields what we were looking for.

\begin{corollary}\label{co:npr}
(i) The consistency of having $n$ measurable cardinals implies the consistency of having $n+1\,$
0-dimensional $T_2$-spaces whose product is irresolvable.

(ii) The consistency of having infinitely many measurable cardinals implies the consistency of having
an $\omega$-sequence $\<X_n : n < \omega\>$ of 0-dimensional $T_2$-spaces such that the product of any
(non-zero) finitely many of them  is irresolvable.
\end{corollary}

By Theorem \ref{tm:econ1} we actually have equiconsistency in (i) for the case $n = 1$, however we do not know this
for $n > 1$. Also, we do not know if equiconsistency holds in (ii).
The simplest remaining problem here is the following.

\begin{problem}\label{pb:3prmes}
Does the existence of an irresolvable product of three (0-dimensional $T_2$) spaces
imply the consistency of having two measurable cardinals?
\end{problem}

\section{On maximal resolvability of products}

Ceder and Pearson proved in \cite{CePe} that if $X$ with $\Delta(X) = \kappa$  is maximally (i.e. $\kappa$)-resolvable then the the product of
$X$ with any space $Y$ of cardinality $\le \kappa^+$ is maximally resolvable. (Of course, this is only of interest if $|Y| = \Delta(Y) = \kappa^+$.)
This prompted them to raise the question if the product of a maximally resolvable space with any space is maximally resolvable.

In \cite{Eck97} Eckertson pointed out that $M(\kappa,\kappa)$ implies a negative answer to the Ceder-Pearson question.
He showed that if $\mathcal{S}$ is a maximal $\kappa$-independent family on $\kappa$ then the the product of
$X(\mathcal{S}) = \<\kappa, \tau(\mathcal{S})\>$ with any space $Y$ of size $< \kappa$ is {\em not} maximally
(i.e. $\kappa$-)resolvable.

By Theorem 1 of \cite{Ku}, $M(\kappa,\kappa)$ implies that $\kappa = |X(\mathcal{S})|$ is strongly inaccessible, hence it is
natural to ask if one could get counterexamples to the Ceder-Pearson question of smaller size. We shall show next that this
is indeed possible by giving two such "small" examples, the first has cardinality $\le 2^{\omega_1}$ and the second has cardinality $\aleph_\omega$.
Just like in Eckertson's case, however, both examples will be obtained from a measurable cardinal. We do not know if this, or actually anything
beyond ZFC, is needed to get a counterexample to the Ceder-Pearson question.

First we give a general result that will be used in getting the first example. We recall from \cite{J} that
$\widehat{c}(X)$ denotes, for any space $X$, the smallest cardinal $\kappa$ such that $X$ does not
contain $\kappa$ pairwise disjoint open sets.

Also, we shall denote the ideal of nowhere dense subsets of $X$ by $\mathcal{N}(X)$ and by $add(\mathcal{N}(X))$
its additivity. The latter is the greatest cardinal such that the union of any fewer nowhere dense sets is nowhere dense.

\begin{theorem}\label{tm:chat}
Let $X$ be an OHI space such that $$\omega < \lambda = \min \big\{\widehat{c}(X), add \big(\mathcal{N}(X) \big)\big\}.$$
Then for any space $Y$ with $|Y| < \lambda$ the product $X \times Y$ is not $\widehat{c}(X)$-resolvable.
\end{theorem}

\begin{proof}
Let $\widehat{c}(X) = \mu$ and consider $X \times Y = \bigcup \{D_\alpha : \alpha < \mu\}$, any $\mu$-partition of the product.
For every point $y \in Y$ and $\alpha < \mu$ we let $E_{y,\alpha} = \{x \in X : \<x, y\> \in D_\alpha\}.$
Clearly, for any $y \in Y$ the sets $\{E_{y,\alpha} : \alpha < \mu \}$ are pairwise disjoint.
Thus, for fixed $y \in Y$, the set $I_y = \{\alpha < \mu : int(E_{y,\alpha}) \ne \emptyset\}$ has cardinality $< \mu$.

But we know that $\widehat{c}(X)$ is always regular, hence $|Y| < \lambda \le \mu$ implies that $I = \bigcup \{I_y : y \in Y\}$
has size $< \mu$ as well, so we may pick an $\alpha \in \mu \setm I$. This means that for each $y \in Y$ the set $E_{y,\alpha}$
has empty interior and so is nowhere dense because $X$ is OHI. But then $|Y| < \lambda \le add \big(\mathcal{N}(X) \big)$
implies that $E = \bigcup \{E_{y,\alpha} : y \in Y\}$ is also nowhere dense, hence there is a non-empty open $U$ in $X$ with $E \cap U = \emptyset$.
However, this clearly implies that $(U \times Y) \cap D_\alpha = \emptyset$, consequently $D_\alpha$ is not dense in $X \times Y$.
\end{proof}

\begin{corollary}\label{co:2adalef1}
Modulo the existence of a measurable cardinal, it is consistent to have
a 0-dimensional $T_2$-space $X$ such that $$\widehat{c}(X) = \omega_2 \le |X| = \Delta(X) \le 2^{\omega_1},$$ and
for any countable space $Y$ the product $X \times Y$ is $\omega_2$-irresolvable. In particular, any such
product is not maximally resolvable, while its factor $Y$ can be chosen to be  maximally resolvable, e.g. $Y = \mathbb{Q}$.
\end{corollary}

\begin{proof}
As we have seen already, we may get $M(\kappa, \omega_1)$ from a measurable cardinal, moreover Kunen proved in \cite{Ku}
that then CH holds and for the smallest $\kappa$ satisfying $M(\kappa, \omega_1)$ we have $\omega_1 < \kappa \le 2^{\omega_1}$.
Also, it follows from CH that if $\mathcal{S}$ is a separating maximal $\omega_1$-independent family on this minimal $\kappa$
then $\widehat{c}(X(\mathcal{S})) = \omega_2$, and $\Delta(X(\mathcal{S})) = |X(\mathcal{S})| = \kappa \ge \omega_2$.
(Actually, it is pointed out in  \cite{Ku} that $\kappa$ is much larger, for instance larger than the first weakly inaccessible cardinal,
but we don't care about that, we only care about the upper bound $2^{\omega_1}$.)

Now, if $X$ is any open OHI subspace of $X(\mathcal{S})$ then we may apply Theorem \ref{tm:chat} to $X$ and $\lambda = \omega_1$
to get the required conclusion. Indeed, for this we only have to check $\omega_1 \le add \big(\mathcal{N}(X) \big)$,
and this trivially follows from the fact that, as we have seen, $X$ has an $\omega_1$-closed base.
\end{proof}

Before giving the second example, we need some preparation that will involve the concept of monotone $\kappa$-resolvability
that has already turned out to be useful in the previous section. It is trivial that $\kappa \in MR(X)$ implies $\kappa \in MR(X \times Y)$
for any space $Y$. Our next result yields a partial converse of this observation.

\begin{theorem}\label{tm:prmr}
Assume that $\kappa = cf(\kappa) \in MR(X \times Y)$, where $|Y| < \kappa$. Then $\kappa \in MR(X)$ as well.
\end{theorem}

\begin{proof}
As was noted above, $\kappa \in MR(X \times Y)$ means that there is a decreasing $\kappa$-sequence $\{D_\alpha : {\alpha}<{\kappa}\}$
of sets {\em dense} in $X \times Y$ with empty intersection. Then the sequence of projections
$\{E_\alpha = \pi_X[D_\alpha] : {\alpha}<{\kappa}\}$ consists of sets dense in $X$ and is also decreasing.
So, it remains to show that $\bigcap \{E_\alpha : {\alpha}<{\kappa}\} = \emptyset$.

To see this, fix any $x \in X$ and note that for any $y \in Y$ there is $\alpha_y < \kappa$ with $\<x, y\> \notin D_{\alpha_y}$.
But $\kappa$ is regular and $|Y| < \kappa$, hence there is an $\alpha < \kappa$ such that $\alpha_y < \alpha$ for all $y \in Y$.
This, however, means that $(\{x\} \times Y) \cap D_\alpha = \emptyset$, hence $x \notin E_\alpha$.
\end{proof}

\begin{corollary}\label{co:alefw}
The consistency of the existence of a measurable cardinal implies that it is consistent to have
a monotonically normal (in short: MN) space $X$ such that $|X| = \Delta(X) = \aleph_\omega$ but
for any countable space $Y$ the product $X \times Y$ is $\omega_1$-irresolvable.
\end{corollary}

\begin{proof}
In \cite{JM}, modulo the existence of a measurable cardinal, a model of ZFC and in it
a MN and (hereditarily) $\omega_1$-irresolvable space $X$ was constructed
such that $|X| = \Delta(X) = \aleph_\omega$. In Theorem 1.5 of \cite{SoSt} it was shown that this space $X$
is not even monotonically  $\omega_1$-resolvable, i.e. $\omega_1 \notin MR(X)$.
Now, it is immediate from Theorem \ref{tm:prmr} that if $Y$ is any countable space then $\omega_1 \notin MR(X \times Y)$
as well, hence $X \times Y$ is $\omega_1$-irresolvable.
\end{proof}

For such a product $X \times Y$ we clearly have $\Delta(X \times Y)  = \aleph_\omega$, hence it
is very far from being maximally resolvable, even if $Y$ is, say $Y = \mathbb{Q}$.

Our next result helps deducing the $\kappa^+$-resolvability of a product from the $\kappa$-resolvability
of one of its factors.

\begin{theorem}\label{tm:kk+}
Assume that $\kappa^+ \in MR(X)$ and $Y$ is $\kappa$-resolvable. Then $X \times Y$ is $\kappa^+$-resolvable.
\end{theorem}

\begin{proof}
Let us start by fixing $\{A_\alpha : {\alpha}<{\kappa^+}\}$, a monotone $\kappa^+$-resolution  of $X$.
As in the proof of Theorem \ref{tm:MRres}, for each point $x \in X$ we let $rk(x) = \min\{\alpha : x \in A_\alpha\}$.

Since $Y$ is $\kappa$-resolvable, we may fix for any ordinal $\alpha < \kappa^+$ a coloring
$c_\alpha : Y \to \alpha$ such that $c_\alpha^{-1}\{\xi\}$ is dense in $Y$ for every $\xi < \alpha$.

Then we define the coloring $d : X \times Y \to \kappa^+$ with the following stipulation:
$$d(x,y) = c_{rk(x)}(y).$$ We claim that  the inverse image $d^{-1}\{\xi\}$ is dense in $X \times Y$
for every $\xi < \kappa^+$.

Indeed, consider any non-empty basic open set $U \times V$ in the product $X \times Y$ and fix $\xi < \kappa^+$.
We may then choose $x \in U$ with $\xi < rk(x)$ and then $y \in V$ such that $c_{rk(x)}(y) = d(x, y) = \xi$,
completing the proof.
\end{proof}

From these two results we easily obtain the following.

\begin{corollary}\label{co:ekv3}
For any cardinal $\kappa$ and space $X$, the following three statements are equivalent.
\begin{enumerate}
  \item $\kappa^+ \in MR(X)$.
  \item $X \times Y$ is $\kappa^+$-resolvable whenever $Y$ is $\kappa$-resolvable.
  \item There is a $\kappa$-resolvable space $Y$ of cardinality $\kappa$ such that $X \times Y$ is $\kappa^+$-resolvable.
\end{enumerate}

\end{corollary}

\begin{proof}
(1) $\Rightarrow$ (2) is just Theorem \ref{tm:kk+}.
(2) $\Rightarrow$ (3) is trivial because there is a $\kappa$-resolvable space $Y$
of cardinality $\kappa$. Finally, (3) $\Rightarrow$ (1) is immediate from Theorem \ref{tm:prmr} applied to $\kappa^+$
because if $X \times Y$ is $\kappa^+$-resolvable then $\kappa^+ \in MR(X \times Y)$.
\end{proof}

From this we get the following improvement on Ceder and Pearson's above cited result.

\begin{corollary}\label{co:CPimpr}
If $X$ is any $\kappa$-resolvable space and $$c(Y) \le \kappa^+ \le \Delta(Y) \le |Y| < \kappa^{\omega},$$
then $X \times Y$ is $\kappa^+$-resolvable.
\end{corollary}

\begin{proof}
As any space is $\kappa^+$-resolvable iff all its neat open subspaces are, it suffices to prove this when $Y$ is neat,
i.e. when $\Delta(Y) = |Y| = \kappa^{+n}$ for some $0 < n < \omega$. Note that in this case we have $\kappa^{+n} \in MR(Y)$.

In Theorem 3.3 of \cite{SoSt} the following "stepping down" result was proved using Ulam matrices:
If $\mu^+ \in MR(Y)$ and $\widehat{c}(Y) \le \mu$ then $\mu \in MR(Y)$ as well. Since $c(Y) \le \kappa^+$, we may apply this $n-1$ times
to conclude that $\kappa^+ \in MR(Y)$. But then $X \times Y$ is $\kappa^+$-resolvable by Corollary \ref{co:ekv3}
(actually with the roles of $X$ and $Y$ switched).
\end{proof}

The most intriguing question concerning the Ceder-Pearson problem that is left open is as follows.

\begin{problem}\label{pb:CP}
Does the existence of a not maximally resolvable product with a maximally resolvable factor
imply the consistency of having a measurable cardinal? Can the existence of such a product
be proven in ZFC?
\end{problem}

We end this section by presenting examples of irresolvable, even submaximal, neat spaces of arbitrarily large cardinality $\kappa$,
whose product with any space of size $\le \kappa$ is maximally (i.e. $\kappa$-)resolvable. In this we shall make use of the
results of \cite{JuSoSz.2005}.

First we need some notation. Given cardinals $\kappa \le \lambda$, we shall denote by $\mathcal{C}(\kappa, \lambda)$ the family
of all $\kappa$-dense subspaces $X$ of the Cantor cube $D(2)^\lambda$ with $|X| = \kappa$. ($X$ is $\kappa$-dense in $D(2)^\lambda$
iff $|X \cap U| \ge \kappa$ for any non-empty open $U$ in $D(2)^\lambda$.)

Now, the following technical lemma is the crucial new ingredient of our promised result.

\begin{lemma}\label{lm:Aay}
Fix an infinite cardinal $\kappa$, moreover let us be given a set $Y$ with $\omega \le |Y| \le \kappa$ and a non-empty family $\mathcal{H} \subs [Y]^\omega$
with $|\mathcal{H}| \le \kappa$ as well.

Then there is a family $$\{A_{\alpha, y} : \<\alpha, y\> \in \kappa \times Y\} \subs \mathcal{C}(\kappa, \kappa)$$
satisfying the following two conditions.
\begin{enumerate}
  \item $\{A_{\alpha, y} : \alpha < \kappa\}$ is pairwise disjoint for every $y \in Y$.

\smallskip

  \item Put $A_{\alpha, H} = \bigcup \{A_{\alpha, y} : y \in H\}$ for any $\<\alpha, H\> \in \kappa \times \mathcal{H}$.\\
  Then for any $0 < n < \omega$ and $\<\eta, h\> \in \kappa^n \times \mathcal{H}^n$ we have
  $$A_{\eta, h} = \bigcap \{A_{\eta(i), h(i)} : i < n\} \in \mathcal{C}(\kappa, \kappa).$$
\end{enumerate}
\end{lemma}

\begin{proof}
Let us start by recalling that the sets $[\varepsilon] = \{x \in D(2)^\kappa : \varepsilon \subs x\}$
with $\varepsilon \in Fn(\kappa, 2)$ form an open base for $D(2)^\kappa$.

Next, fix a linear ordering $\prec$ of $Y$. Then, for any $0 < n < \omega$ and $j \in [Y]^n$, for every $i < n$ we let $j(i)$
denote the $i$th member of $j$ in its $\prec$-increasing order.

Let us now consider  the set of triples $$\mathbb{T} = \bigcup \{\kappa^n \times [Y]^n \times Fn(\kappa, 2) : 0 < n < \omega\}.$$
Clearly, we have $|\mathbb{T}| = \kappa$, so we may fix an enumeration $$\mathbb{T} = \{\<\eta_\zeta, j_\zeta, \varepsilon_\zeta\> : \zeta < \kappa\}$$
such that $|\{\zeta : \<\eta, j, \varepsilon\> = \<\eta_\zeta, j_\zeta, \varepsilon_\zeta\>\}| = \kappa$ for every $\<\eta, j, \varepsilon\> \in \mathbb{T}$.
Note that for every $\zeta < \kappa$ we have $n_\zeta \in \omega \setm \{0\}$ such that $\<\eta_\zeta, j_\zeta\> \in \kappa^{n_\zeta} \times [Y]^{n_\zeta}$.

After this we may use simple transfinite recursion to define points $x_\zeta \in D(2)^\kappa$ such that
$x_\zeta \in [\varepsilon_\zeta] \setm \{x_\xi : \xi < \zeta\}$ for every $\zeta < \kappa$.
Finally, for every $\<\alpha, y\> \in \kappa \times Y$ we set
$$A_{\alpha, y} = \{x_\zeta : \exists\, i < n_\zeta \text{ such that } \eta_\zeta(i) = \alpha \text{ and } j_\zeta(i) = y\}.$$
It is obvious from our construction that then $|[\varepsilon] \cap A_{\alpha, y}| = \kappa$ for every $\varepsilon \in Fn(\kappa, 2)$,
hence $A_{\alpha, y} \in \mathcal{C}(\kappa, \kappa)$.

Now, assume that $x_\zeta \in A_{\alpha, y}$. Then, by definition, there is a, clearly unique, $i < n_\zeta$
such that $ j_\zeta(i) = y$, hence we also must have $\eta_\zeta(i) = \alpha$. In other words, no $x_\zeta$
may belong to $A_{\alpha, y}  \cap A_{\beta, y}$ if $\alpha \ne \beta$, hence we have verified (1).

Finally, to check (2), consider any $\<\eta, h\> \in \kappa^n \times \mathcal{H}^n$ for some $0 < n < \omega$.
We may then find pairwise distinct $y_i \in h(i)$ for each $i < n$, hence $j = \{y_i : i < n\} \in [Y]^n$.
Thus for every $\varepsilon \in Fn(\kappa, 2)$ the triple $\<\eta, j, \varepsilon \> \in \mathbb{T}$,
while it is clear that $x_\zeta \in A_{\eta, h} \cap [\varepsilon]$ whenever $\<\eta, j, \varepsilon\> = \<\eta_\zeta, j_\zeta, \varepsilon_\zeta\>$.
Consequently we do have $A_{\eta, h} \in \mathcal{C}(\kappa, \kappa).$
\end{proof}

We are now ready to formulate and prove the promised result.

\begin{theorem}\label{tm:submax}
For every cardinal $\kappa = \kappa^\omega$ there is a {\em submaximal} space $X \in \mathcal{C}(\kappa, 2^\kappa)$
such that the product of $X$ with any space of cardinality $\le \kappa$ is $\kappa$-resolvable, i.e. maximally resolvable.
\end{theorem}

\begin{proof}
Since $\kappa = \kappa^\omega$, we may apply Lemma \ref{lm:Aay} with $Y = \kappa$ and $\mathcal{H} = [\kappa]^\omega$ to obtain
$\{A_{\alpha, y} : \<\alpha, y\> \in \kappa \times Y\} \subs \mathcal{C}(\kappa, \kappa)$ satisfying conditions (1) and (2).
(We keep denoting the second indexes by $y$ for better readability.)
Let $E$ be any dense subset of $D(2)^{2^\kappa \setm \kappa}$ of cardinality $\kappa$ and for each $\<\alpha, y\> \in \kappa \times Y$
set $E_{\alpha, y} = A_{\alpha, y} \times E.$ Clearly, each $E_{\alpha, y} \in \mathcal{C}(\kappa, 2^\kappa)$.
By (2) for any pair $\<\eta, h\> \in \kappa^n \times \mathcal{H}^n$ with $0 < n < \omega$ we also have
$E_{\eta, h} = A_{\eta, h} \times E \in \mathcal{C}(\kappa, 2^\kappa)$.

By the main Theorem 3.3 of \cite{JuSoSz.2005} then there is a bijection $f$ from $S = \bigcup \{E_{\alpha, y} : \<\alpha, y\> \in \kappa \times Y\}$
onto some NODEC $X \in \mathcal{C}(\kappa, 2^\kappa)$ such that each $D_{\eta, h} = f[E_{\eta, h}] \in \mathcal{C}(\kappa, 2^\kappa)$,
moreover the topology of $X$ is $\mathcal{D}$-forced, where $$\mathcal{D} = \{D_{\eta, h} : \<\eta, h\> \in \cup\{\kappa^n \times \mathcal{H}^n : 0 < n < \omega\}.$$
This means that for every set $A$ that is dense in some non-empty open $U$ in $X$ there is a non-empty open $V \subs U$ and a
$D_{\eta, h} \in \mathcal{D}$ such that $V \cap D_{\eta, h} \subs A$, see Fact 2.4 in \cite{JuSoSz.2005}.

To see that $X$ is OHI, consider any non-empty open $U$ in $X$ and two dense subsets, say $A$ and $B$, of $U$.
Then there are some $\<\eta, h\> \in \kappa^n \times \mathcal{H}^n$ and a non-empty open $V \subs U$ such that $V \cap D_{\eta, h} \subs A$.
But $B \cap V$ is dense in $V$, hence
there are $\<\vartheta, k\> \in \kappa^m \times \mathcal{H}^m$ and a non-empty open $W \subs V$ such that $W \cap D_{\vartheta, k} \subs B$.
Now, using that $$D_{\eta, h} \cap D_{\vartheta, k} = D_{\eta \smallfrown \vartheta, h \smallfrown k} \in \mathcal{D}$$
is also dense in $X$, we conclude that $\emptyset \ne W \cap D_{\eta \smallfrown \vartheta, h \smallfrown k} \subs A \cap B$.
So, any two dense subsets of $U$ meet, i.e. $U$ is irresolvable. Since a space is submaximal iff it is both NODEC and OHI,
$X$ is indeed submaximal.

Note that by definition, we have $$X = f[S] = \cup \{D_{\alpha, y} = f[E_{\alpha, y}] : \<\alpha, y\> \in \kappa \times Y\}.$$
For every $\alpha < \kappa$ let us now define $$D_\alpha = \cup \{D_{\alpha, y} \times \{y\} : y \in Y\} \subs X \times Y.$$
It is then immediate from (1) that $\alpha \ne \beta$ implies $D_\alpha \cap D_\beta = \emptyset$.

Now, for any $\alpha < \kappa$ and $H \in \mathcal{H} = [Y]^\omega$ we have $D_{\alpha, H} = \cup \{D_{\alpha, y} : y \in H\} \in \mathcal{D}$,
hence $D_{\alpha, H}$ is dense in $X$. This implies that for every $\alpha < \kappa$ we have $D_\alpha \cap (U \times H) \ne \emptyset$
whenever $U$ is non-empty open in $X$ and $H \in [Y]^\omega$.

Now, if $Z$ is any space with $|Z| \le \kappa$ then we may assume that its underlying set is included in $Y$, and as every non-empty open set in $Z$
is infinite, it clearly follows that $\{D_\alpha \cap (X \times Z) : \alpha < \kappa\}$ are pairwise disjoint dense sets in $X \times Z$.
Thus $X \times Z$ is indeed $\kappa$-resolvable.
\end{proof}

We do not know if Theorem \ref{tm:submax} remains valid without the assumption $\kappa = \kappa^\omega$. However, we do have a "local"
version of it which does not require this assumption. To formulate it, we introduce a new cardinal function $d_\omega$ that we
call $\omega$-density. We say that $\mathcal{H} \subs [Y]^\omega$ is dense in a space $Y$ if for every non-empty open set $U$ in $Y$
there is $H \in \mathcal{H}$ such that $H \subs U$.

\begin{definition}
For any space $X$ we let $$d_\omega(X) = \min \{|\mathcal{H}| : \mathcal{H} \text{ is dense in }X\}.$$
\end{definition}

Note that by our convention $d_\omega(X)$ is well-defined for every space $X$.
We clearly have $d(X) \le d_\omega(X) \le \pi(X)$ for any $X$. Also, if $X$ is $T_1$ then $d_\omega(X) \le d(X)^\omega$.

\begin{theorem}\label{tm:submaxQ}
Let $Y$ be any space.Then
for every cardinal $\kappa \ge \max\{|Y|, d_\omega(Y)\}$ there is a submaximal space $X \in \mathcal{C}(\kappa, 2^\kappa)$
such that the product  $X \times Y$ is $\kappa$-resolvable, i.e. maximally resolvable.
\end{theorem}

Since the proof of this theorem, based on Lemma \ref{lm:Aay}, is completely analogous to the proof of Theorem \ref{tm:submax},
we leave it to the reader.

Note that for the space $\mathbb{Q}$ of the rationals $\max\{|\mathbb{Q}|, d_\omega(\mathbb{Q})\} = \omega$,
hence for every $\kappa \ge \omega$ we have a submaximal $X \in \mathcal{C}(\kappa, 2^\kappa)$ with $X \times \mathbb{Q}$
maximally resolvable. This is in contrast with Corollaries \ref{co:2adalef1} and \ref{co:alefw} which provide products
with $\mathbb{Q}$ that are not maximally resolvable.

\bibliographystyle{plain}

\begin{thebibliography}{10}

    \bibitem{CePe}
    J.~Ceder and T.~Pearson.
    \newblock On products of maximally resolvable spaces.
    \newblock {\em Pacific J. Math.}, 22:31--45, 1967.

    \bibitem{Eck97}
    Frederick~W. Eckertson.
    \newblock Resolvable, not maximally resolvable spaces.
    \newblock {\em Topology Appl.}, 79(1):1--11, 1997.

    \bibitem{je74}
    Ronald~Bj\"{o}rn Jensen.
    \newblock Measurable cardinals and the {${\rm GCH}$}.
    \newblock In {\em Axiomatic set theory ({P}roc. {S}ympos. {P}ure {M}ath.,
      {V}ol. {XIII}, {P}art {II}, {U}niv. {C}alifornia, {L}os {A}ngeles, {C}alif.,
      1967)}, pages 175--178, 1974.

    \bibitem{J}
    Istv\'{a}n Juh\'{a}sz.
    \newblock {\em Cardinal functions in topology---ten years later}, volume 123 of
      {\em Mathematical Centre Tracts}.
    \newblock Mathematisch Centrum, Amsterdam, second edition, 1980.

    \bibitem{JM}
    Istvan Juh\'{a}sz and Menachem Magidor.
    \newblock On the maximal resolvability of monotonically normal spaces.
    \newblock {\em Israel J. Math.}, 192(2):637--666, 2012.

    \bibitem{JuSoSz.2005}
    Istv\'{a}n Juh\'{a}sz, Lajos Soukup, and Zolt\'{a}n Szentmikl\'{o}ssy.
    \newblock {$\mbb D$}-forced spaces: a new approach to resolvability.
    \newblock {\em Topology Appl.}, 153(11):1800--1824, 2006.

    \bibitem{Ku}
    Kenneth Kunen.
    \newblock Maximal {$\sigma $}-independent families.
    \newblock {\em Fund. Math.}, 117(1):75--80, 1983.

    \bibitem{Kb}
    Kenneth Kunen.
    \newblock {\em Set theory}, volume~34 of {\em Studies in Logic (London)}.
    \newblock College Publications, London, 2011.

    \bibitem{KST}
    Kenneth Kunen, Andrzej Szyma\'{n}ski, and Franklin Tall.
    \newblock Baire irresolvable spaces and ideal theory.
    \newblock {\em Ann. Math. Sil.}, (14):98--107, 1986.

    \bibitem{KT}
    Kenneth Kunen and Franklin Tall.
    \newblock On the consistency of the non-existence of baire irresolvable spaces.
    \newblock Manuscript, Topology Atlas,\url{http://at.yorku.ca/v/a/a/a/27.htm}.

    \bibitem{Ma73}
    V.~I. Malyhin.
    \newblock Products of ultrafilters, and indecomposable spaces.
    \newblock {\em Mat. Sb. (N.S.)}, 90(132):106--116, 166, 1973.

    \bibitem{SoSt}
    Lajos Soukup and Adrienne Stanley.
    \newblock Resolvability in c.c.c. generic extensions.
    \newblock {\em Comment. Math. Univ. Carolin.}, 58(4):519--529, 2017.

    \end{thebibliography}

\end{document}